\documentclass[12pt]{article}

\usepackage{amsmath, epsfig, cite}
\usepackage{amssymb}
\usepackage{amsfonts}
\usepackage{latexsym}
\usepackage{graphicx}
\usepackage{url}
\usepackage{color}
\usepackage{cite}

\newtheorem{thm}{Theorem}[section]

\newtheorem{cor}[thm]{Corollary}
\newtheorem{defi}[thm]{Definition}
\newtheorem{lem}[thm]{Lemma}
\newtheorem{conj}[thm]{Conjecture}

\newcommand{\qed}{{\hfill\rule{4pt}{7pt}}}
\def\pf{\noindent {\it Proof.} }

\numberwithin{equation}{section}

\makeatletter \@addtoreset{equation}{section} \makeatother

\begin{document}
\rule{0cm}{1cm}
\begin{center}
{\large  On the Maximum Number of $k$-Hooks of Partitions of $n$ }
\end{center}
 \vskip 2mm \centerline{Anna R.B. Fan$^1$, Harold R.L. Yang$^2$ and Rebecca T.
 Yu$^3$ }

\begin{center}
Center for Combinatorics, LPMC-TJKLC\\
Nankai University\\ Tianjin 300071
P.R.China\\

\vskip 2mm
$^1$fanruice@mail.nankai.edu.cn, $^2$yangruilong@mail.nankai.edu.cn,
$^3$yuting\_shuxue@mail.nankai.edu.cn
\end{center}

{\noindent \bf Abstract.}
Let $\alpha_k(\lambda)$ denote the number of   $k$-hooks  in a   partition $\lambda$  and let
 $b(n,k)$ be the maximum  value of $\alpha_k(\lambda)$ among  partitions of $n$.
 Amdeberhan
 posed a conjecture on the generating function
  of $b(n,1)$. We give a proof of this conjecture. In general,
  we obtain a formula that can be used to determine $b(n,k)$. This leads to a
 generating function formula for $b(n,k)$.
We introduce the notion of
 nearly $k$-triangular partitions.  We show that for any $n$,
 there is a nearly $k$-triangular partition which can be transformed into
 a partition of $n$ that attains the maximum number
  of $k$-hooks. The operations for the transformation enable us to compute the
  number $b(n,k)$.
\vskip 3mm



\section{Introduction}

The objective of this paper is to derive a generating function formula for the maximum number of $k$-hooks
 in the Young diagrams of partitions of $n$. For $k=1$, the problem was
 posed by Amdeberhan \cite{Tewodros}. Let $\alpha_1(\lambda)$ be the number of $1$-hooks in the partition $\lambda$, or equivalently, the number of distinct
  parts in $\lambda$. 
Let
\[b_n=\max\{\alpha_1(\lambda) \colon \lambda \in P(n)\},\]
where $P(n)$ denotes the set of  partitions of $n$.

 Amdeberhan \cite{Tewodros} posed the following conjecture.
 \begin{conj}\label{conj}
 We have
\begin{align}\label{bngenerate}
\sum_{n\geq 0} b_n q^n=\frac{1}{1-q}\left(\frac{(q^2;q^2)^2_{\infty}}{(q;q)_{\infty}}-1\right),
\end{align}
where $(q,q)_\infty=(1-q)(1-q^2)(1-q^3)\cdots$.
\end{conj}

Following the notation $\alpha_1(\lambda)$ of
  Amdeberhan, we use $\alpha_k(\lambda)$ to denote  the number of $k$-hooks in $\lambda$, and let
\[b(n,k)=\max\{\alpha_k(\lambda)\colon \lambda \in P(n)\}.\] The main result of this
paper is the following generating function formula for $b(n,k)$.

\begin{thm}\label{thm3}
For $n\geq 0$ and $k\geq 1$, we have
\begin{align}
\sum_{n\geq 0}b(n,k) q^n=\frac{1}{1-q}\left(\sum_{t\geq 1}q^{{t \choose 2}k^2}\frac{1-q^{tk^2}}{1-q^{tk}}-1\right).\label{generatebnk}
\end{align}
\end{thm}

Clearly, Theorem \ref{thm3} reduces to Theorem \ref{bngenerate} when $k=1$.
Pak \cite{Pak}  gave the following generating function of  the statistic $\alpha_1(\lambda)$,
where he
used  $\gamma(\lambda)$ to denote $\alpha_1(\lambda)$:
\begin{equation}
\sum_{n\geq 0}\sum_{\lambda \in P(n)}\gamma(\lambda)q^n=\frac{q}{(1-q)(q;q)_\infty}.\label{pakeq}
\end{equation}

In general, the statistic $\alpha_k(\lambda)$ has been
studied by Han \cite[Eq. $1.5$]{Han}. More precisely,  he derived the following generating
function formula:
\begin{equation}\label{haneq}
\sum_{n\geq 0}\sum_{\lambda \in  P(n) }x^{\alpha_k(\lambda)} q^{|\lambda|}=\prod_{j\geq 1}\frac{(1+(x-1)q^{kj})^k}
{1-q^j}.
\end{equation}

Taking   logarithms of both sides of \eqref{haneq} and differentiating  with respect to $x$, we
obtain the following relation by setting $x=1$:
\begin{equation}\label{han2}
\sum_{n\geq 0}\sum_{\lambda\in P(n)}\alpha_k(\lambda)q^{n}=\frac{kq^k}{(1-q^k)(q;q)_{\infty}}.
\end{equation}
It can be seen that identity \eqref{han2} becomes \eqref{pakeq} when $k=1$.

 To prove Theorem \ref{thm3}, we introduce a class of partitions, called nearly $k$-triangular partitions.
\begin{defi}\label{tmk}
For fixed $m\geq 0$ and $k\geq 1$,   let $m=lk+r$, where $0\leq r\leq k-1$.
Let  $T_m^{(k)}$ denote the nearly $k$-triangular partition with $m$
parts as given by
$$T_m^{(k)}=({\mathop{ \underbrace{(l+1)k,\ldots,(l+1)k}} \limits_{r},\mathop{\underbrace{lk,\ldots,lk}}\limits_{k}
,\ldots,\mathop{\underbrace{2k,\ldots,2k}}\limits_{k},\mathop{\underbrace{k,\ldots,k}}\limits_{k}}).$$
\end{defi}

\begin{figure}[h]
\begin{center}
\setlength{\unitlength}{1mm}
\begin{picture}(50,58)
\multiput(0,54)(0,-6){9}{\line(1,0){6}}
\multiput(0,54)(0,-6){8}{\line(0,-1){6}}
\multiput(6,54)(0,-6){8}{\line(0,-1){6}}
\multiput(6,54)(0,-6){9}{\line(1,0){6}}
\multiput(12,54)(0,-6){8}{\line(0,-1){6}}
\multiput(12,54)(0,-6){9}{\line(1,0){6}}
\multiput(18,54)(0,-6){8}{\line(0,-1){6}}
\multiput(18,54)(0,-6){6}{\line(1,0){6}}
\multiput(24,54)(0,-6){5}{\line(0,-1){6}}
\multiput(24,54)(0,-6){6}{\line(1,0){6}}
\multiput(30,54)(0,-6){5}{\line(0,-1){6}}
\multiput(30,54)(0,-6){6}{\line(1,0){6}}
\multiput(36,54)(0,-6){5}{\line(0,-1){6}}
\multiput(36,54)(0,-6){3}{\line(1,0){6}}
\multiput(42,54)(0,-6){2}{\line(0,-1){6}}
\multiput(42,54)(0,-6){3}{\line(1,0){6}}
\multiput(48,54)(0,-6){2}{\line(0,-1){6}}
\multiput(48,54)(0,-6){3}{\line(1,0){6}}
\multiput(54,54)(0,-6){2}{\line(0,-1){6}}
\multiput(44,48.5)(-6,-6){8}{*}
\thicklines\put(0,24){\line(1,0){18}}\put(0,24){\line(0,-1){18}}
           \put(18,24){\line(0,-1){18}}\put(0,6){\line(1,0){18}}
           \put(18,42){\line(1,0){18}}\put(18,42){\line(0,-1){18}}
           \put(36,42){\line(0,-1){18}}\put(18,24){\line(1,0){18}}
           \put(36,54){\line(1,0){18}}\put(36,54){\line(0,-1){12}}
           \put(54,54){\line(0,-1){12}}\put(36,42){\line(1,0){18}}

\end{picture}
\caption{A nearly $3$-triangular partition $T_8^{(3)}$.}\label{figure}
\end{center}
\end{figure}
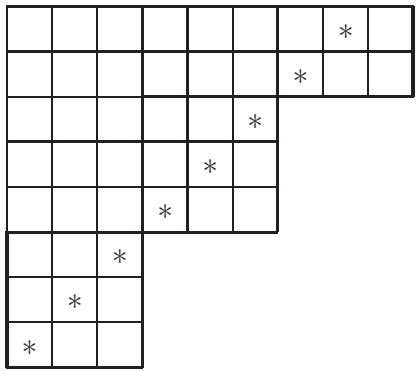

For example, Figure \ref{figure} illustrates a nearly $3$-triangular
 partition with eight parts. Throughout the paper, we use the symbol $*$ to mark  these cells with hook
 length $k$.

It can be seen that there is exactly one cell in each row of
 $T_m^{(k)}$ with hook length $k$.

Let us recall some basic  notation and terminology on partitions as used in
\cite{George}. {\it A partition} $\lambda$ of a positive integer
$n$ is a finite nonincreasing sequence of positive integers
$\lambda=(\lambda_1,\,\lambda_2,\ldots,\,\lambda_r)$
 such that
$\sum_{i=1}^r\lambda_i=n.$  The entries $\lambda_i$ are called the parts
of $\lambda$. The number of parts of  $\lambda$ is called the
length of $\lambda$, denoted by $l(\lambda).$ The weight of
$\lambda$ is the sum of parts, denoted by $|\lambda|.$

\begin{figure}[h]
\begin{center}
\setlength{\unitlength}{1mm}
\begin{picture}(160,24)
\put(40,8){\hfill\rule{11.5pt}{11.5pt}}
\put(44,8){\hfill\rule{11.5pt}{11.5pt}}
\put(40,4){\hfill\rule{11.5pt}{11.5pt}}
\multiput(36,20)(0,-4){6}{\line(1,0){4}}
\multiput(36,20)(0,-4){5}{\line(0,-1){4}}
\multiput(40,20)(0,-4){5}{\line(0,-1){4}}
\multiput(40,20)(0,-4){5}{\line(1,0){4}}
\multiput(44,20)(0,-4){4}{\line(0,-1){4}}
\multiput(44,20)(0,-4){4}{\line(1,0){4}}
\multiput(48,20)(0,-4){3}{\line(0,-1){4}}
\multiput(48,20)(0,-4){2}{\line(1,0){4}}
\multiput(52,20)(0,-4){1}{\line(0,-1){4}}
\multiput(52,20)(0,-4){2}{\line(1,0){4}}
\multiput(56,20)(0,-4){1}{\line(0,-1){4}}
\multiput(56,20)(0,-4){2}{\line(1,0){4}}
\multiput(60,20)(0,-4){1}{\line(0,-1){4}}
\multiput(92,20)(0,-4){6}{\line(1,0){4}}
\multiput(92,20)(0,-4){5}{\line(0,-1){4}}
\multiput(96,20)(0,-4){5}{\line(0,-1){4}}
\multiput(96,20)(0,-4){5}{\line(1,0){4}}
\multiput(100,20)(0,-4){4}{\line(0,-1){4}}
\multiput(100,20)(0,-4){4}{\line(1,0){4}}
\multiput(104,20)(0,-4){3}{\line(0,-1){4}}
\multiput(104,20)(0,-4){2}{\line(1,0){4}}
\multiput(108,20)(0,-4){1}{\line(0,-1){4}}
\multiput(108,20)(0,-4){2}{\line(1,0){4}}
\multiput(112,20)(0,-4){1}{\line(0,-1){4}}
\multiput(112,20)(0,-4){2}{\line(1,0){4}}
\multiput(116,20)(0,-4){1}{\line(0,-1){4}}
\put(97.1,7.5){*}
\put(105.1,15.6){*}
\put(93.1,3.5){*}
\end{picture}
\caption{A $3$-hook  and the three cells with hook length 3.}\label{figure1.1}
\end{center}
\end{figure}
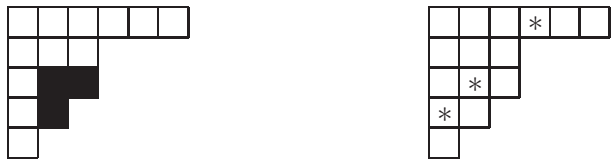

A partition can be represented by a Young diagram. For each
 cell $u$ in  the Young diagram of  $\lambda$,  we define
 the hook length $h_u(\lambda)$ of $u$  by
  the number of cells $v$ in the Young diagram of $\lambda$ such that $v = u$,
   or $v$ appears below $u$ in the same column, or $v$ lies to the right of $u$ in the same row.
    A hook with length $k$ is called a $k$-hook, see Figure \ref{figure1.1}.

This paper is organized as follows. In Section $2$, we give a
formula that can be used to determine $b_n$  in Theorem \ref{lembn}.
We then derive the generating function
for $b_n$.  In Section $3$, we define two operations $D_i$ and $P_i$ on Young diagrams. Using these
operations one can transform a partition
  $\lambda$ with $m$ $k$-hooks
  into a nearly $k$-triangular partition $T_s^{(k)}$, where $s\geq m$. This leads to a proof of Theorem \ref{thm1}.
In Section 4, we define an operation $Q_j$ on  Young diagrams.
We obtain a formula for  $b(n,k)$ as well as a formula for the generating function.

\section{Proof of Conjecture \ref{conj}}

In this section, we give a proof of Conjecture 1.1.

 \begin{thm}\label{lembn}  Assume that $m$ is a nonnegative integer and $n$ is an integer such that
${m+1\choose 2}\leq n \leq {m+2\choose 2}-1.$ Then we have $b_n=m$.
\end{thm}

\pf Recall that  $\alpha_1(\lambda)$   is
 the number of distinct parts of $\lambda $ and $b_n$ is the maximum value of $\alpha_1(\lambda)$ when $\lambda$
  ranges over partitions of $n$.
 We claim that we have $b_n<m+1$ for any $n< {m+2 \choose 2}$.  Assume that $\lambda$  is a partition
 with $m+1$ distinct parts.
It is clear that
 \[|\lambda| \geq 1+2+\cdots +(m+1)={m+2\choose 2}.\]
In other words,  if $n<{m+2\choose 2}$ then we have
$b_n<m+1$. So the claim is verified.

Next we proceed to show that $b_n\geq m$ for any $n\geq{m+1\choose 2}$.
 Let
\[\lambda=(m,m-1,\cdots,2,1^{n-{m+1\choose 2}+1}).\]
Clearly, $\lambda$ has  $m$ distinct parts and $|\lambda|=n$. Thus $b_n \geq m$ if $n \geq {m+1\choose 2}$.
 So we reach the conclusion
 that $b_n=m$ for ${m+1\choose 2}\leq n \leq {m+2\choose 2}-1.$  This completes the proof. \qed

We are ready to prove Conjecture \ref{conj} with the aid of Theorem  \ref{lembn}.

\noindent {\it{ Proof of Conjecture \ref{conj}.}} First, we may express the
generating function of $b_n$ in terms of the generating function of $b_{n+1}-b_n$.
More precisely, we have
\begin{align}
(1-q)\sum_{n\geq 0}b_nq^n=\sum_{n\geq 0}(b_{n+1}-b_n)q^{n+1}. \label{difbn}
\end{align}
To compute $b_{n+1}-b_n$, we denote the interval
  $[{m+1\choose 2}, {m+2\choose 2}-1]$ by $I_m$. By Theorem \ref{lembn},
we see that $b_n$ is determined by the interval which $n$ lies in.
There are two cases:

\noindent Case 1: $n$ and $n+1$ lie in the same interval $I_m$.
 Then we have $b_{n+1}=b_n=m$, that is, $b_{n+1}-b_n=0.$

\noindent Case 2: $n$ and $n+1$ lie in two consecutive intervals $I_{m}$ and $I_{m+1}$,
that is, $n={m+2 \choose 2}-1$. In this case, we have $b_n=m$ and $b_{n+1}=m+1$.
So we have  $b_{n+1}-b_n=1.$

Combining the above two cases, we find that
\[
  \sum_{n\geq 0}(b_{n+1}-b_n)q^{n+1}=\sum_{m\geq 0}q^{m+2\choose 2}
= \sum_{m\geq 0}q^{m+1\choose 2}-1.
\]
By Gauss' identity \cite[Eq. $1.4.10$]{George2}
\begin{align}
\sum_{n=0}^{\infty}q^{ {n+1 \choose 2}}=\frac{(q^2;q^2)_{\infty}}{(q;q^2)_{\infty}},\label{psi}
\end{align}
we obtain that
\begin{equation} \label{temp} (1-q)\sum_{n\geq 0}b_nq^n= \sum_{m\geq 0}q^{m+1\choose 2}-1=\frac{(q^2;q^2)^2_{\infty}}{(q;q)_{\infty}}-1.
\end{equation}
Thus, identity \eqref{bngenerate} can be deduced  from \eqref{temp} by dividing both sides
by $(1-q)$. This completes the proof.
\qed

\section{Nearly $k$-triangular partitions}

In this section, we introduce the structure of nearly $k$-triangular partitions, and we show that
such a partition has minimum weight given the number of $k$-hooks. This property
will be used in the next section to determine $ b(n,k)$.

For $m\geq 0$ and $k\geq 1$, the weight of the nearly $k$-triangular partition $T_m^{(k)}$
is given by
\begin{equation}\label{defi-NM}
t(m,k)= m\left(\left\lfloor \frac{m}{k}\right\rfloor +1\right)k-{\lfloor\frac{m}{k}\rfloor +1\choose 2}k^2.
  \end{equation} It can be seen that there is exactly one cell in each row of $T_m^{(k)}$
  that has hook length $k$. Hence $T_m^{(k)}$ is a partition with $m$ $k$-hooks. The following
  theorem states that  $T_m^{(k)}$ has minimum weight among partitions with $m$ $k$-hooks.

\begin{thm}\label{thm1}
 For $m\geq 0$ and $k\geq 1$, if  $\lambda$ is a partition with $m$ $k$-hooks, then we have
$|\lambda|\geq  t(m,k).$
\end{thm}

To prove Theorem \ref{thm1}, we introduce two operations  $D_i$ and $P_i$ defined on Young diagrams.
They can also be considered as operations on partitions.
 We shall  show that one can transform  a partition $\lambda$ with $m$ $k$-hooks into
 a nearly $k$-triangular partition  $T_m^{(k)}$ by applying the operations $D_i$ and $P_i$.

The   operations $D_i$ and $P_i$ are defined as follows.
Let $\lambda=(\lambda_1,\lambda_2,\ldots, \lambda_r)$ be a partition.
The operation $D_i$ means to remove the $i$-th row of the
Young diagram of $\lambda$.
The operation $P_i$ applies to partitions $\lambda$ for which  $\lambda_i>\lambda_{i+1}$.
More precisely, $P_i(\lambda)$  is  obtained from $\lambda$ via the following steps.
Assume that the cells with hook length $k$ are  marked by $*$.

\noindent Step 1.  Remove the last cell $u$  from the $i$-th
row of the Young diagram of $\lambda$, and denote the resulting partition by $\mu$.
If the Young diagram of $\lambda$ contains no marked cell in the column occupied by $u$,
  then we set $P_i(\lambda)=\mu$;

\noindent Step 2.  At this step,  there is a cell of hook length $k$ in
the column of $\lambda$ that contains $u$.
Denote this marked cell by $w_1$ and assume that it is in the $j$-th row of $\lambda$.
Evidently, the marked cell $w_1$ in $\lambda$ is of hook length $k-1$ in $\mu$.
There are two cases:

\noindent Case 1: $\mu_{j}=\mu_{j-1}$.
Notice that the cell $w_1'$ above $w_1$ is of hook length $k$ in $\mu$.
That is to say that $w_1'$ is a marked cell in $\mu$.
We set $P_i(\lambda)=\mu$;

\noindent Case 2:  $\mu_{j}<\mu_{j-1}$.
We add one cell $v$ at the end of the $j$-th row in $\mu$ and denote the resulting partition by $\nu$.
Clearly,  $w_1$ is also a marked cell in $\nu$.
If the Young diagram of $\mu$ contains no marked cell in the $(\mu_j+1)$-th column,
  then we set $P_i(\lambda)=\nu$. Otherwise,  go to the next step.

\noindent Step 3. There is a marked cell  in the $(\mu_j+1)$-th column
 of $\mu$. Let $w_2$ denote this marked cell and suppose that it is  in the $h$-th row.
Evidently,  the cell $w_2$ has hook length $k+1$ in $\nu$.
There are two cases:

\noindent Case 1:  $\nu_{k}=\nu_{k+1}$.
Notice that the cell $w_2'$  below $w_2$ is of hook length $k$ in $\nu$.
In this case, we set $P_i(\lambda)=\nu$;

\noindent Case 2: $\nu_{k}>\nu_{k+1}$.
We remove the last cell $u'$ in  the $h$-th row of $\nu$ and denote the resulting partition by $\mu'$.
Now, $w_2$ is also a marked cell in $\mu'$.
If the Young diagram of $\nu$ contains no marked cell in the column occupied by $u'$,
we set $P_i(\lambda)=\mu'$. Otherwise,
there is a marked cell  in the  column occupied by $u'$ in $\nu$.
We set $u=u'$, $\lambda=\nu$, $\mu=\mu'$
and go back to Step 2.

\begin{figure}[h]
\begin{center}
\setlength{\unitlength}{1.5mm}
\begin{picture}(40,42)
\put(0,38){\line(0,-1){38}}
\put(0,38){\line(1,0){48}}
\put(0,0){\line(1,0){6}}
\put(6,0){\line(0,1){4}}
\put(6,4){\line(1,0){10}}
\put(16,4){\line(0,1){4}}
\put(0,8){\line(1,0){22}}
\put(0,10){\line(1,0){22}}
\put(20,8){\line(0,1){8}}
\put(22,8){\line(0,1){8}}
\put(20.25,8.5){\scriptsize{$-$}}\put(23,8){$u$}\put(-3,8.3){\small{$i$}}
\put(0,14){\line(1,0){28}}
\put(0,16){\line(1,0){28}}
\put(26,14){\line(0,1){10}}
\put(28,14){\line(0,1){10}}
\put(20.25,14.25){$*$}\put(20,17){\small{$w_1$}}
\put(26.25,14.5){\scriptsize{$+$}}\put(29,14){\small{$v$}}\put(-3,14.3){\small{$j$}}
\put(0,22){\line(1,0){32}}
\put(0,24){\line(1,0){32}}
\put(30,22){\line(0,1){8}}
\put(32,22){\line(0,1){8}}
\put(26.25,22.25){$*$}\put(26,25){\small{$w_2$}}
\put(30.25,22.5){\scriptsize{$-$}}\put(33,22){\small{$u'$}}\put(-3,22.3){\small{$h$}}
\put(0,28){\line(1,0){36}}
\put(0,30){\line(1,0){36}}
\put(36,28){\line(0,1){7}}
\put(36,35){\line(1,0){12}}
\put(48,35){\line(0,1){3}}
\put(30.25,28.25){$*$}
\end{picture}
\end{center}
\caption{The operation $P_i$.}\label{general}
\end{figure}
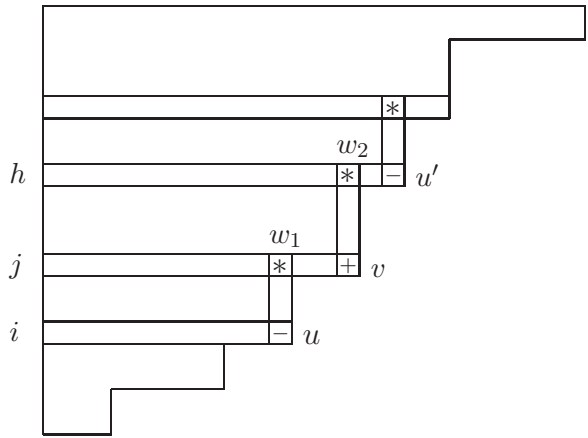

Figure \ref{general} gives an  illustration of the operation $P_i$. The cells with the symbol $-$
are the removed cells and the cells with the symbol $+$ are the added cells.

The following property of the operation $P_i$ is easy to verify, and hence the proof is omitted.
\begin{lem}\label{lemma1}
Let $\lambda=(\lambda_1,\lambda_2,\ldots,\lambda_r)$
be a partition such that  $\lambda_i>\lambda_{i+1}$ and let $\alpha_k(\lambda,i)$
denote the number of marked cells in the $i$-th row of $\lambda$. Then we have
$|\lambda|  \geq |P_i(\lambda)|$ and
\begin{align*}
\alpha_k(\lambda) - \alpha_k(\lambda,i) \leq \alpha_k(P_i(\lambda))- \alpha_k(P_i(\lambda),i).
\end{align*}
\end{lem}

We are now in a position to present a proof of  Theorem  \ref{thm1}
by using the operations $P_i$ and $D_i$.

\noindent {\it{Proof of Theorem \ref{thm1}.}}
Let $\lambda=(\lambda_1,\lambda_2,\ldots,\lambda_r)$ be a partition with $m$ $k$-hooks.
We shall give a  procedure to transform $\lambda$ into $T_s^{(k)}$ where $s \geq m$ by using the operations $D_i$ and $P_i$.
Notice that $D_i$ decreases the weight of a partition and $P_i$ either keeps the weight or
decreases the weight by one. Moreover, it can be seen that during the process of the
transformation $D_i$ preserves the number of $k$-hooks and $P_i$ keeps or increases the
number of $k$-hooks.
Hence we arrive at the conclusion that $|\lambda| \geq t(s,k) \geq t(m,k)$.

Clearly, we have $r\geq m$.
We aim to construct a sequence $\beta^{(r)}, \beta^{(r-1)},\ldots, \beta^{(1)}, \beta^{(0)}$ of nearly
$k$-triangular partitions starting with $T_r^{(k)}$ and
ending with $T_{s}^{(k)}$, where $s\geq m$.
In the construction of $T_{s}^{(k)}$ from $\lambda$, we denote the intermediate partitions by $\lambda^{(r)},\lambda^{(r-1)},$ $\ldots,\lambda^{(1)},\lambda^{(0)}$ with $\lambda^{(r)}=
\lambda$ and  $\lambda^{(0)}=T_s^{(k)}$. We compare the partitions $\beta^{(i)}$  and $\lambda^{(i)}$ to construct $\beta^{(i-1)}$ and $\lambda^{(i-1)}$ by the following process, where,
 $\beta^{(i)}$ and $\lambda^{(i)}$ have  the same number of parts and  $\lambda^{(i)}$ and $\beta^{(i)}$ differ
 only in the first $i$ rows.

First, we compare the last entry of $\lambda^{(r)}$ with the last entry of $\beta^{(r)}$. Recall that
 $\lambda^{(r)}=\lambda$ and $\beta^{(r)}=T_r^{(k)}$. There are two cases:

\noindent Case 1: The $r$-th entry of $\lambda^{(r)}$ is equal to the $r$-th entry of $\beta^{(r)}$.
Then we set $\beta^{(r-1)}=T_{r}^{(k)}$ and $\lambda^{(r-1)}=\lambda^{(r)}$.

\noindent Case 2: The $r$-th entry of $\lambda^{(r)}$ is not equal to the $r$-th entry of $\beta^{(r)}$.
There are two subcases:

\noindent Case 2.1:  $ \lambda^{(r)}_r<\beta^{(r)}_r=k$. Clearly, there are no cells
of hook length $k$ in the $r$-th row of the Young diagram of $\lambda^{(r)}$.
Applying the operation $P_r$ to $\lambda^{(r)}$, the $r$-th entry of $\lambda^{(r)}$ decreases by one.
 In view of Lemma \ref{lemma1}, we have $\alpha_k(P_r(\lambda^{(r)}))\geq \alpha_k(\lambda)$ and $|P_r(\lambda^{(r)})|\leq |\lambda|$. Hence, applying the operation $P_r$ to $\lambda^{(r)}$ iteratively,
 more precisely, $\lambda_r^{(r)}$ times, we obtain a new partition $\mu$ with $r-1$ parts. It is clear that $\alpha_k(\mu)\geq \alpha_k(\lambda^{(r)})$ and $|\mu|\leq |\lambda^{(r)}|$. We set $\beta^{(r-1)}=T_{r-1}^{(k)}$ and $\lambda^{(r-1)}=\mu$.

\noindent Case 2.2: $\lambda^{(r)}_r>\beta^{(r)}_r=k$.
Evidently, there is a cell with hook length $k$ in the
$r$-th row of $\lambda^{(r)}$.
Applying the operation $P_r$ to $\lambda^{(r)}$  $(\lambda_r^{(r)}-k)$ times,
we obtain a new partition $\mu$. It is clear that the $r$-th entry of $\mu$ equals the $r$-th entry of $\beta^{(r)}$.
By Lemma \ref{lemma1},
we find that $\alpha_k(\mu)\geq \alpha_k(\lambda^{(r)})$ and $|\mu|\leq |\lambda^{(r)}|$.
We set $\beta^{(r-1)}=T_{r}^{(k)}$ and $\lambda^{(r-1)}=\mu$.

Then we proceed to compare the $i$-th entry of $\lambda^{(i)}$ with the $i$-th entry of $\beta^{(i)}$ for $r-1
\geq i \geq 1$,
where we assume that $\lambda^{(i)}$ and $\beta^{(i)}$ have been constructed by the above procedure.
Assume that $\lambda^{(i)}$ has $t$ parts and $\beta^{(i)}=T_t^{(k)}$,
There are two cases:

\noindent Case 1: The $i$-th entry of $\lambda^{(i)}$ is equal to the $i$-th entry of $\beta^{(i)}$.
Then we set $\beta^{(i-1)}=T_t^{(k)}$ and $\lambda^{(i-1)}=\lambda^{(i)}$.

\noindent Case 2: The $i$-th entry of $\lambda^{(i)}$ is not equal to the $i$-th entry of $\beta^{(i)}$.
 There are three subcases:

\noindent Case 2.1:  $i{\equiv} t \pmod{k}$ and $ \lambda^{(i)}_i-\lambda^{(i)}_{i+1}<k$. In this case,
 let $d=\lambda^{(i)}_i-\lambda^{(i)}_{i+1}$. Clearly, there are no cells
of hook length $k$ in the $i$-th row of the Young diagram of $\lambda^{(i)}$,
see Figure \ref{case2}. Applying the operation $P_i$ to $\lambda^{(i)}$  $d$ times,
we obtain a new partition $\mu$ with $\mu_i=\mu_{i+1}$. It is clear that
there are no marked cells in the $i$-th row of $\mu$.
In view of Lemma \ref{lemma1},
we see that  $\alpha_k(\mu)\geq \alpha_k(\lambda^{(i)})$ and $|\mu|\leq |\lambda^{(i)}|$.

Now we apply the operation $D_i$ to $\mu$ to obtain a new partition $\nu$ with $t-1$ parts.
Since there are no marked cells in the area $A$ of $\mu$, namely,  $\{(p,q):1\leq p\leq i-1, 1\leq q\leq \mu_i\}$, see Figure \ref{case2}. The positions of the marked cells of $\mu$
stay unchanged in $\nu$ with respect to the operation $D_i$.
It follows that $\alpha_k(\nu)= \alpha_k(\mu)$ and $|\nu|< |\mu|$.
 This implies that $\alpha_k(\nu)\geq \alpha_k(\lambda^{(i)})$ and $|\nu|< |\lambda^{(i)}|$.
Notice that the partitions $\nu$ and $T_{t-1}^{(k)}$ differ only in the first $i-1$ rows.
We set $\beta^{(i-1)}=T_{t-1}^{(k)}$ and $\lambda^{(i-1)}=\nu$.

\begin{figure}[h]
\begin{center}
\setlength{\unitlength}{1mm}
\begin{picture}(20,36)
\put(-12,0){\line(1,0){12}}
\put(-12,0){\line(0,1){32}}{\thicklines\put(-12,16){\line(0,1){16}}}
                           {\thicklines\put(-12,16){\line(1,0){12}}}
\put(0,0){\line(0,1){12}}
\put(-12,32){\line(1,0){40}}{\thicklines\put(-12,32){\line(1,0){12}}}
                            {\thicklines\put(0,32){\line(0,-1){16}}}
\put(20,16){\line(0,1){4}}
\put(20,20){\line(1,0){8}}
\put(28,20){\line(0,1){12}}

\put(-12,12){\line(1,0){20}}
\put(-12,16){\line(1,0){32}}

\multiput(0,12)(0,1){20}{$\line(0,1){0.5}$}
\multiput(-12,12)(1,0){24}{$\line(1,0){0.5}$}
\multiput(-12,12)(0,1){20}{$\line(0,1){0.5}$}
\multiput(-12,24)(1,0){36}{$\line(1,0){0.5}$}
\multiput(12,12)(0,1){20}{$\line(0,1){0.5}$}
\multiput(-12,32)(1,0){36}{$\line(1,0){0.5}$}
\multiput(24,24)(0,1){8}{$\line(0,1){0.5}$}

\multiput(-11.2,0.7)(3.8,3.8){3}{$*$}
\multiput(-12,12)(2,0){9}{$\line(1,1){4}$}

\put(0,12){\line(0,1){20}}
\put(-8,22.5){$A$}

\put(8,12){\line(0,1){4}}

\put(9.5,12.5){\small the $i$-th row}

\end{picture}
\end{center}
\caption{The case for $i{\equiv} t \pmod{k}$ and $ \lambda^{(i)}_i-\lambda^{(i)}_{i+1}<k$.}\label{case2}
\end{figure}

\noindent Case 2.2: $i{\equiv} t \pmod{k}$ and $\lambda^{(i)}_i-\lambda^{(i)}_{i+1}>k$.
Denote $d=\lambda^{(i)}_i-\lambda^{(i)}_{i+1}-k$.  Evidently, there is a cell with hook length $k$ in the
$i$-th row of $\lambda^{(i)}$, see Figure
\ref{case3}. Applying the operation $P_i$ to $\lambda^{(i)}$  $d$ times,
we obtain a new partition $\mu$.
It is clear that there is also a marked cell in the $i$-th row.
In view of Lemma \ref{lemma1}, we see that $\alpha_k(\mu)\geq \alpha_k(\lambda^{(i)})$ and $|\mu|\leq |\lambda^{(i)}|$.
Now the partitions $\mu$ and $T_{t}^{(k)}$ differ only in the first $i-1$ rows.
We set $\beta^{(i-1)}=T_{t}^{(k)}$ and $\lambda^{(i-1)}=\mu$.

\begin{figure}[h]
\begin{center}
\setlength{\unitlength}{1mm}
\begin{picture}(20,36)
\put(-12,0){\line(1,0){12}}
\put(-12,0){\line(0,1){32}}
\put(0,0){\line(0,1){12}}
\put(-12,32){\line(1,0){40}}
\put(20,16){\line(0,1){4}}
\put(20,20){\line(1,0){8}}
\put(28,20){\line(0,1){12}}
\put(-12,16){\line(1,0){32}}

\put(-12,12){\line(1,0){28}}
\put(16,12){\line(0,1){4}}

\multiput(0,12)(0,1){20}{$\line(0,1){0.5}$}
\multiput(-12,12)(1,0){24}{$\line(1,0){0.5}$}
\multiput(-12,12)(0,1){20}{$\line(0,1){0.5}$}
\multiput(-12,24)(1,0){36}{$\line(1,0){0.5}$}
\multiput(12,12)(0,1){20}{$\line(0,1){0.5}$}
\multiput(-12,32)(1,0){36}{$\line(1,0){0.5}$}
\multiput(24,24)(0,1){8}{$\line(0,1){0.5}$}

\multiput(-11.2,0.7)(3.8,3.8){3}{$*$}
\multiput(-12,12)(2,0){13}{$\line(1,1){4}$}

\put(6,12.5){$*$}

\put(17,12.5){\small the $i$-th row}

\end{picture}
\end{center}
\caption{ The case for $i{\equiv} t \pmod{k}$ and $\lambda^{(i)}_i-\lambda^{(i)}_{i+1}>k$.}\label{case3}
\end{figure}

\begin{figure}[h]
\begin{center}
\setlength{\unitlength}{1mm}
\begin{picture}(20,36)
\put(-12,0){\line(1,0){12}}
\put(-12,0){\line(0,1){32}}
\put(0,0){\line(0,1){12}}
\put(-12,32){\line(1,0){40}}
\put(20,16){\line(0,1){4}}
\put(20,20){\line(1,0){8}}
\put(28,20){\line(0,1){12}}
\put(-12,16){\line(1,0){32}}

\put(-12,12){\line(1,0){24}}
\put(12,12){\line(0,1){4}}
\put(-12,20){\line(1,0){32}}
\put(1,12.5){$*$}

\multiput(0,12)(0,1){20}{$\line(0,1){0.5}$}
\multiput(-12,12)(1,0){24}{$\line(1,0){0.5}$}
\multiput(-12,12)(0,1){20}{$\line(0,1){0.5}$}
\multiput(-12,24)(1,0){36}{$\line(1,0){0.5}$}
\multiput(12,12)(0,1){20}{$\line(0,1){0.5}$}
\multiput(-12,32)(1,0){36}{$\line(1,0){0.5}$}
\multiput(24,24)(0,1){8}{$\line(0,1){0.5}$}

\multiput(-11.2,0.7)(3.8,3.8){3}{$*$}
\multiput(-12,16)(2,0){15}{$\line(1,1){4}$}

\put(21,16.5){\small the $i$-th row}

\end{picture}
\end{center}
\caption{The case for $i\not\equiv t \pmod{k}$.}\label{case4}
\end{figure}
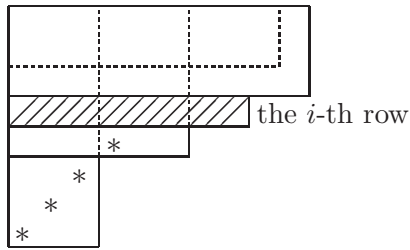

\noindent Case 2.3: $i\not\equiv t \pmod{k}$. In this case, we have
 $\lambda^{(i)}_i-\lambda^{(i)}_{i+1}>0$,
see Figure \ref{case4}. Let $d=\lambda^{(i)}_i-\lambda^{(i)}_{i+1}$. Applying the operation  $P_i$
to $\lambda^{(i)}$ $d$ times, we obtain a new partition  $\mu$.
Clearly,  there is a marked cell in the $i$-th row of $\mu$.
By  Lemma \ref{lemma1}, we deduce that $\alpha_k(\mu)\geq \alpha_k(\lambda^{(i)})$ and $|\mu|\leq |\lambda^{(i)}|$.
Now, the partitions $\mu$ and $T_{t}^{(k)}$ differ only in the first $i-1$ rows.
We set $\beta^{(i-1)}=T_t^{(k)}$ and $\lambda^{(i-1)}=\mu$.

Repeat the above process, we eventually obtain a nearly $k$-triangular partition $\lambda^{(0)}=\beta^{(0)}=T_s^{(k)}$.
In the construction of $\lambda^{(0)}$ from $\lambda^{(r)}$, the operations $P_i$ and $D_i$ imply
the following relations:
\[m=\alpha_k(\lambda^{(r)})\leq \cdots \leq\alpha_k(\lambda^{(i)})\leq \alpha_k(\lambda^{(i-1)})\leq\cdots\leq\alpha_k(\lambda^{(0)})=s\]
and
\[|\lambda|=|\lambda^{(r)}|\geq \cdots \geq|\lambda^{(i)}|\geq |\lambda^{(i-1)}|\geq\cdots\geq|\lambda^{(0)}|=t(s,k).\]
Since $s\geq m$, we have $|\lambda|\geq |\lambda^{(0)}|= t(s,k)\geq t(m,k)$. This completes the proof.  \qed

As a consequence of Theorem \ref{thm1},  we obtain the following upper bound on $b(n,k)$.

\begin{cor}\label{cor-r} Assume  $m\geq 0$ and $k>0$ and $n$ is a nonnegative integer  such that
 $ n <t(m+1,k) $,
  then $ b(n,k)\leq m.$
\end{cor}

\begin{figure}[h]
\begin{center}
\setlength{\unitlength}{1mm}
\begin{picture}(50,80)
\multiput(-32,76)(0,-4){8}{\line(1,0){4}}
\multiput(-32,76)(0,-4){7}{\line(0,-1){4}}
\multiput(-28,76)(0,-4){7}{\line(0,-1){4}}
\multiput(-28,76)(0,-4){8}{\line(1,0){4}}
\multiput(-24,76)(0,-4){7}{\line(0,-1){4}}
\multiput(-24,76)(0,-4){8}{\line(1,0){4}}
\multiput(-20,76)(0,-4){7}{\line(0,-1){4}}
\multiput(-20,76)(0,-4){4}{\line(1,0){4}}
\multiput(-16,76)(0,-4){3}{\line(0,-1){4}}
\multiput(-16,76)(0,-4){3}{\line(1,0){4}}
\multiput(-12,76)(0,-4){2}{\line(0,-1){4}}
\multiput(-12,76)(0,-4){3}{\line(1,0){4}}
\multiput(-8,76)(0,-4){2}{\line(0,-1){4}}
\multiput(-8,76)(0,-4){3}{\line(1,0){4}}
\multiput(-4,76)(0,-4){2}{\line(0,-1){4}}
\multiput(-4,76)(0,-4){2}{\line(1,0){4}}
\multiput(0,76)(0,-4){1}{\line(0,-1){4}}
\multiput(0,76)(0,-4){2}{\line(1,0){4}}
\multiput(4,76)(0,-4){1}{\line(0,-1){4}}
\multiput(4,76)(0,-4){2}{\line(1,0){4}}
\multiput(8,76)(0,-4){1}{\line(0,-1){4}}
\multiput(-31,47.5)(4,4){3}{*}
\put(-15,67.5){*}
\put(-3,71.5){*}
\put(-15,57){\scriptsize$|\lambda^{(4)}|=33$}
\put(-15,53){\scriptsize$l(\lambda^{(4)})=7$}
\put(-15,49.5){\scriptsize$\alpha_3(\lambda^{(4)})=5$}

\multiput(4.2,76)(0,-1){5}{$\line(0,-1){0.4}$}
\multiput(-8,71.8)(1,0){13}{$\line(1,0){0.4}$}
\multiput(-7.8,72)(0,-1){12}{$\line(0,-1){0.4}$}
\multiput(-20,59.8)(1,0){13}{$\line(1,0){0.4}$}
\multiput(-19.8,60)(0,-1){12}{$\line(0,-1){0.4}$}
\multiput(-32,47.8)(1,0){12}{$\line(1,0){0.4}$}

\put(10,62){\vector(1,0){18}}
\put(16,64){\small$D_4$}
\put(12,58){\small Case 2.1}
\multiput(44,76)(0,-4){7}{\line(1,0){4}}
\multiput(44,76)(0,-4){6}{\line(0,-1){4}}
\multiput(48,76)(0,-4){6}{\line(0,-1){4}}
\multiput(48,76)(0,-4){7}{\line(1,0){4}}
\multiput(52,76)(0,-4){6}{\line(0,-1){4}}
\multiput(52,76)(0,-4){7}{\line(1,0){4}}
\multiput(56,76)(0,-4){6}{\line(0,-1){4}}
\multiput(56,76)(0,-4){4}{\line(1,0){4}}
\multiput(60,76)(0,-4){3}{\line(0,-1){4}}
\multiput(60,76)(0,-4){3}{\line(1,0){4}}
\multiput(64,76)(0,-4){2}{\line(0,-1){4}}
\multiput(64,76)(0,-4){3}{\line(1,0){4}}
\multiput(68,76)(0,-4){2}{\line(0,-1){4}}
\multiput(68,76)(0,-4){3}{\line(1,0){4}}
\multiput(72,76)(0,-4){2}{\line(0,-1){4}}
\multiput(72,76)(0,-4){2}{\line(1,0){4}}
\multiput(76,76)(0,-4){1}{\line(0,-1){4}}
\multiput(76,76)(0,-4){2}{\line(1,0){4}}
\multiput(80,76)(0,-4){1}{\line(0,-1){4}}
\multiput(80,76)(0,-4){2}{\line(1,0){4}}
\multiput(84,76)(0,-4){1}{\line(0,-1){4}}
\multiput(45,51.5)(4,4){3}{*}
\put(61,67.5){*}
\put(73,71.5){*}
\put(80,57){\scriptsize$|\lambda^{(3)}|=30$}
\put(80,53){\scriptsize$l(\lambda^{(3)})=6$}
\put(80,49.5){\scriptsize$\alpha_3(\lambda^{(3)})=5$}

\multiput(68.2,76)(0,-1){12}{$\line(0,-1){0.4}$}
\multiput(56,63.8)(1,0){13}{$\line(1,0){0.4}$}
\multiput(56.2,64)(0,-1){12}{$\line(0,-1){0.4}$}
\multiput(44.2,51.8)(1,0){12}{$\line(1,0){0.4}$}

\put(70,50){\vector(0,-1){16}}
\put(71,43){\small$P_3$}
\put(71,38){\small$D_3$}
\put(55,40){\small Case 2.1}
\multiput(60,26)(0,-4){6}{\line(1,0){4}}
\multiput(60,26)(0,-4){5}{\line(0,-1){4}}
\multiput(64,26)(0,-4){5}{\line(0,-1){4}}
\multiput(64,26)(0,-4){6}{\line(1,0){4}}
\multiput(68,26)(0,-4){5}{\line(0,-1){4}}
\multiput(68,26)(0,-4){6}{\line(1,0){4}}
\multiput(72,26)(0,-4){5}{\line(0,-1){4}}
\multiput(72,26)(0,-4){3}{\line(1,0){4}}
\multiput(76,26)(0,-4){2}{\line(0,-1){4}}
\multiput(76,26)(0,-4){3}{\line(1,0){4}}
\multiput(80,26)(0,-4){2}{\line(0,-1){4}}
\multiput(80,26)(0,-4){3}{\line(1,0){4}}
\multiput(84,26)(0,-4){2}{\line(0,-1){4}}
\multiput(84,26)(0,-4){3}{\line(1,0){4}}
\multiput(88,26)(0,-4){2}{\line(0,-1){4}}
\multiput(88,26)(0,-4){2}{\line(1,0){4}}
\multiput(92,26)(0,-4){1}{\line(0,-1){4}}
\multiput(92,26)(0,-4){2}{\line(1,0){4}}
\multiput(96,26)(0,-4){1}{\line(0,-1){4}}
\multiput(96,26)(0,-4){2}{\line(1,0){4}}
\multiput(100,26)(0,-4){1}{\line(0,-1){4}}
\multiput(61,5.5)(4,4){3}{*}
\put(77,17.5){*}
\put(89,21.5){*}
\put(77,3){\scriptsize$|\lambda^{(2)}|=26$}
\put(77,-1){\scriptsize$l(\lambda^{(2)})=5$}
\put(77,-4.5){\scriptsize$\alpha_3(\lambda^{(2)})=5$}

\multiput(84.2,26)(0,-1){9}{$\line(0,-1){0.4}$}
\multiput(72,17.8)(1,0){13}{$\line(1,0){0.4}$}
\multiput(72.2,18)(0,-1){12}{$\line(0,-1){0.4}$}
\multiput(60,5.8)(1,0){12}{$\line(1,0){0.4}$}

\put(55,15){\vector(-1,0){16}}
\put(45,17){\small$P_2$}
\put(40,11){\small Case 2.2}
\multiput(2,26)(0,-4){6}{\line(1,0){4}}
\multiput(2,26)(0,-4){5}{\line(0,-1){4}}
\multiput(6,26)(0,-4){5}{\line(0,-1){4}}
\multiput(6,26)(0,-4){6}{\line(1,0){4}}
\multiput(10,26)(0,-4){5}{\line(0,-1){4}}
\multiput(10,26)(0,-4){6}{\line(1,0){4}}
\multiput(14,26)(0,-4){5}{\line(0,-1){4}}
\multiput(14,26)(0,-4){3}{\line(1,0){4}}
\multiput(18,26)(0,-4){2}{\line(0,-1){4}}
\multiput(18,26)(0,-4){3}{\line(1,0){4}}
\multiput(22,26)(0,-4){2}{\line(0,-1){4}}
\multiput(22,26)(0,-4){3}{\line(1,0){4}}
\multiput(26,26)(0,-4){2}{\line(0,-1){4}}
\multiput(26,26)(0,-4){2}{\line(1,0){4}}
\multiput(30,26)(0,-4){1}{\line(0,-1){4}}
\multiput(30,26)(0,-4){2}{\line(1,0){4}}
\multiput(34,26)(0,-4){1}{\line(0,-1){4}}
\multiput(34,26)(0,-4){2}{\line(1,0){4}}
\multiput(38,26)(0,-4){1}{\line(0,-1){4}}
\multiput(38,26)(0,-4){2}{\line(1,0){4}}
\multiput(42,26)(0,-4){1}{\line(0,-1){4}}
\multiput(3,5.5)(4,4){3}{*}
\put(15,17.5){*}
\put(31,21.5){*}
\put(17,3){\scriptsize$|\lambda^{(1)}|=25$}
\put(17,-1){\scriptsize$l(\lambda^{(1)})=5$}
\put(17,-4.5){\scriptsize$\alpha_3(\lambda^{(1)})=5$}

\multiput(26.2,26)(0,-1){9}{$\line(0,-1){0.4}$}
\multiput(14,17.8)(1,0){13}{$\line(1,0){0.4}$}
\multiput(14.2,18)(0,-1){12}{$\line(0,-1){0.4}$}
\multiput(2,5.8)(1,0){12}{$\line(1,0){0.4}$}

\put(-7,15){\vector(-1,0){16}}
\put(-19,17){\small$(P_1)^4$}
\put(-22,11){\small Case 2.3}
\multiput(-48,26)(0,-4){6}{\line(1,0){4}}
\multiput(-48,26)(0,-4){5}{\line(0,-1){4}}
\multiput(-44,26)(0,-4){5}{\line(0,-1){4}}
\multiput(-44,26)(0,-4){6}{\line(1,0){4}}
\multiput(-40,26)(0,-4){5}{\line(0,-1){4}}
\multiput(-40,26)(0,-4){6}{\line(1,0){4}}
\multiput(-36,26)(0,-4){5}{\line(0,-1){4}}
\multiput(-36,26)(0,-4){3}{\line(1,0){4}}
\multiput(-32,26)(0,-4){2}{\line(0,-1){4}}
\multiput(-32,26)(0,-4){3}{\line(1,0){4}}
\multiput(-28,26)(0,-4){2}{\line(0,-1){4}}
\multiput(-28,26)(0,-4){3}{\line(1,0){4}}
\multiput(-24,26)(0,-4){2}{\line(0,-1){4}}
\multiput(-31,21.5)(-4,-4){5}{*}
\put(-33,3){\scriptsize$|\lambda^{(0)}|=21$}
\put(-33,-1){\scriptsize$l(\lambda^{(0)})=5$}
\put(-33,-4.5){\scriptsize$\alpha_3(\lambda^{(0)})=5$}

\multiput(-23.8,26)(0,-1){9}{$\line(0,-1){0.4}$}
\multiput(-36,17.8)(1,0){13}{$\line(1,0){0.4}$}
\multiput(-35.8,18)(0,-1){12}{$\line(0,-1){0.4}$}
\multiput(-48,5.8)(1,0){12}{$\line(1,0){0.4}$}

\end{picture}
\vskip 10mm
\caption{ $\lambda=(10,7,4,3,3,3,3)$ and  $T_5^{(3)}=(6,6,3,3,3).$
}\label{transform}
\end{center}
\end{figure}
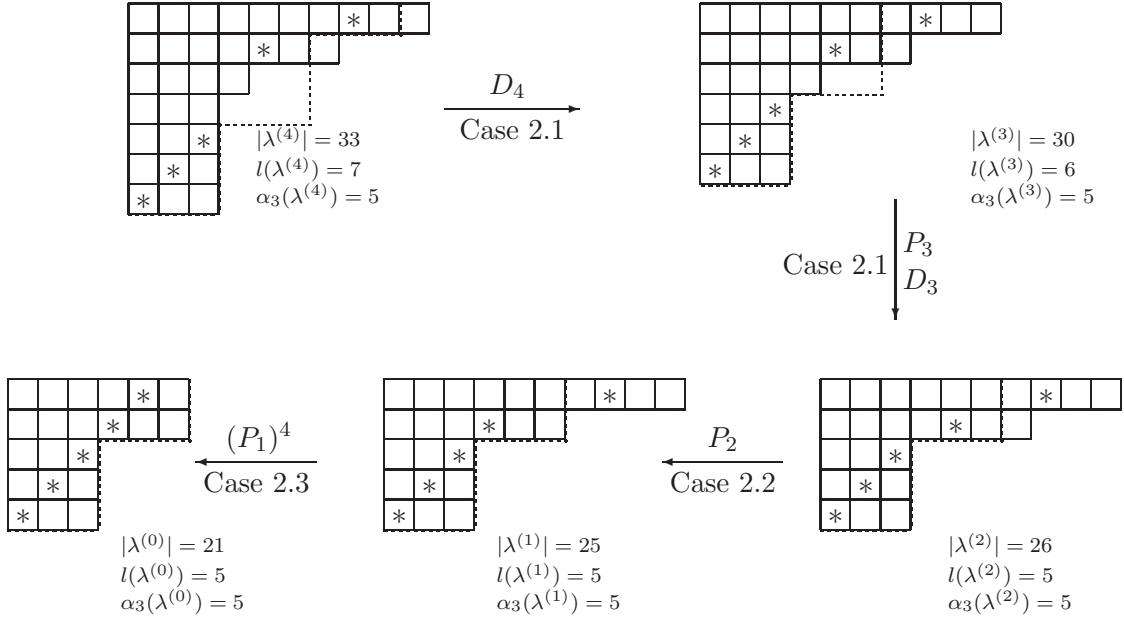

Figure
\ref{transform}   illustrates the transformation from
$\lambda=(10,7,4,3,3,3,3)$ to a nearly 3-triangular partition $T_5^{(3)}=(6,6,3,3,3)$.
It can be checked that  both $\lambda$ and $T_5^{(3)}$
 have five $3$-hooks and $|\lambda|>t(5,3)=21.$

\section{Proof of Theorem \ref{thm3}}

In this section, we show that the number $b(n,k)$, that is, the maximum number
of $k$-hooks among partitions of $n$, can be determined by the number $t(m,k)$, namely,
the weight of the nearly $k$-triangular partition $T_m^{(k)}$. In the previous section,
we have obtained an upper bound on $b(n,k)$. To determine $b(n,k)$, we shall give a lower
bound on $b(n,k)$, as stated in Theorem \ref{thm2}.

 \begin{thm}\label{cor}Assume $m\geq 0 $, $k\geq 1$ and $n$ is a nonnegative integer  such that
$t(m,k)\leq n \leq t(m+1,k)-1$, then we have $b(n,k)=m.$
\end{thm}

Note that Theorem \ref{cor} reduces to Theorem \ref{lembn} when $k=1$.
 The following theorem gives a lower bound on $b(n,k)$.

\begin{thm}\label{thm2}  Assume  $m\geq 0$ and $k>0$ and $n$ is a nonnegative integer  such that
 $n \geq t(m,k)$,
  then we have $b(n,k)\geq m.$
\end{thm}

To prove Theorem \ref{thm2}, we  introduce an operation
$Q_j$ defined on Young diagrams. In fact,  $Q_j$ is same as
the operation $P_i$ except for the first step.
 Let $\lambda=(\lambda_1,\lambda_2,\ldots, \lambda_r)$ be a partition. The operation $Q_j$ applies to partitions $\lambda$ for which  $\lambda_{j-1}>\lambda_j$.
More precisely, $Q_j(\lambda)$ is constructed via the following steps.

\noindent Step 1.
 Add a cell $v$
at the end of the $j$-th row of the Young diagram of $\lambda$, and denote the resulting partition by $\mu$.
If the Young diagram of $\lambda$ contains no marked cells in the $(\lambda_j+1)$-th column of $\lambda$,
  then we set $Q_j(\lambda)=\mu$;

\noindent Step 2. At this step,  it is clear that there is only one cell of hook length $k$ in
the $(\lambda_j+1)$-th column of $\lambda$.
Denote this marked cell by $w$ and assume that it is in the $h$-th row of $\lambda$.
Note that the marked cell $w$ in $\lambda$ is of hook length $k+1$ in $\mu$.
There are two cases:

\noindent Case 1:  $\mu_h=\mu_{h+1}$.
Clearly, the cell $w'$  below $w$ is of hook length $h$ in $\mu$.
In this case, we set $Q_j(\lambda)=\mu$

\noindent Case 2:  $\mu_h>\mu_{h+1}$.
We apply the operation $P_h$ to $\mu$ and denote the resulting partition by $\nu$.
 It is easily seen that the $h$-th entry of $\mu$ decreases by one. Consequently, $w$ has  hook length $k$ in $\nu$. We set $Q_j(\lambda)=\nu$.

For a partition $\lambda=(\lambda_1,\lambda_2,\ldots,\lambda_r)$, we regard the  $(r+1)$-th entry as $0$ when we apply $Q_{r+1}$ to $\lambda$. Under this convention,
 $Q_{r+1}$ increases the number of parts of $\lambda$ by one.

The following property of the operation $Q_j$ is similar to Lemma \ref{lemma1}. The proof is omitted.
\begin{lem}\label{lemma2}
Let $\lambda=(\lambda_1,\lambda_2,\ldots,\lambda_r)$ be a partition such that $\lambda_j<\lambda_{j-1}$.
Then we have
\begin{equation}\label{alphak}
\alpha_k(\lambda) - \alpha_k(\lambda,j) \leq \alpha_k(Q_j(\lambda))- \alpha_k(Q_j(\lambda),j)
\end{equation}
and
\begin{equation}\label{weight}
 |\lambda|\leq |Q_j(\lambda)|\leq |\lambda|+1.
 \end{equation}

\end{lem}

We are now ready to prove Theorem \ref{thm2} by using the operation $Q_j$.

\noindent{\it{Proof of Theorem \ref{thm2}.}} For each $n\geq t(m,k)$,
it suffices to show that there exists a partition $\lambda$ of $n$
with at least $m$ $k$-hooks.  We proceed by  induction on $n$.

First, when $n=t(m,k)$, the nearly $k$-triangular partition $T_m^{(k)}$ is
a partition of $t(m,k)$ with $m$ $k$-hooks, that is the theorem holds for $n=t(m,k)$.

We now assume that there is a partition $\lambda=(\lambda_1,\lambda_2,\ldots,\lambda_r)$ of $N$ with $s$ $k$-hooks,
where $N\geq t(m,k)$ and $s\geq m$. The following procedure gives the
construction of a partition  of $N+1$
with at least $s$ $k$-hooks.

 Apply $Q_{r+1}$ to $\lambda$ and denote the resulting partition by $\mu^{(1)}$.
 From Lemma \ref{lemma2}, we see that $\alpha_k(\mu^{(1)})\geq s$ and $|\lambda|\leq |\mu^{(1)}|\leq|\lambda|+1.$
There are two cases:

\noindent Case 1: $|\mu^{(1)}|=|\lambda|+1=N+1$. Then $\mu^{(1)}$ is a
 partition of $N+1$ with at least $s$ $k$-hooks.

\noindent Case 2: $|\mu^{(1)}|=|\lambda|=N$. We continue to
construct a sequence of partitions,
$$\mu^{(2)}=Q_{r+2}(\mu^{(1)}),
\ldots,\mu^{(i+1)}=Q_{r+i+1}(\mu^{(i)}),\ldots ,\mu^{(k-1)}=Q_{r+k-1}(\mu^{(k-2)}).$$
It is clear that $\mu^{(i)}$ has $r+i$ parts and contains at least $i$ parts equal to one.
 It follows from  \eqref{alphak} that
 $\alpha_k(\mu^{(i)})\geq s$ for $2 \leq i \leq k-1$. There are two subcases:

\noindent Case 2.1: There exists a partition $\mu^{(i)}$ such that $|\mu^{(i)}|=N+1$, where
 $2 \leq i \leq k-1$.
Then $\mu^{(i)}$ is a partition of $N+1$ with at least $s$ $k$-hooks.

\noindent Case 2.2: $|\mu^{(i)}|=N$ for $2 \leq i \leq k-1$.
 Now, we construct $\mu^{(k)}$ from $\mu^{(k-1)}$. Recall that $\mu^{(k-1)}$ contains
  at least $k-1$ parts equal to
 one with at least $s$ $k$-hooks. Set $\mu^{(k)}$ to be the partition
 obtained from $\mu^{(k-1)}$ by adding one as a new part.
 Now, we have $|\mu^{(k)}|=N+1$.
 Moreover, there are at least $k$ parts equal one and
 there is a cell of hook length $k$ in the first column of $\mu^{(k)}$.
 It can be seen that the positions of the marked cells in the other column of $\mu^{(k-1)}$ stay unchanged in $\mu^{(k)}$. This yields  $\alpha_k(\mu^{(k)})\geq \alpha_k(\mu^{(k-1)})\geq s$.
 Thus we have obtained a partition  $\mu^{(k)}$ of $N+1$ with at least $s$ $k$-hooks.\qed

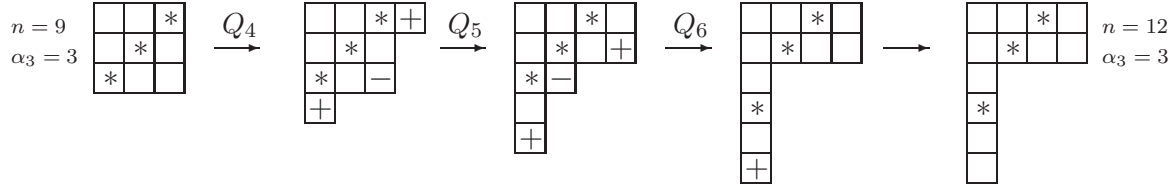
\begin{figure}[h]
\begin{center}
\setlength{\unitlength}{1mm}
\begin{picture}(26,24)
\multiput(-52,20)(0,-4){4}{\line(1,0){4}}
\multiput(-52,20)(0,-4){3}{\line(0,-1){4}}
\multiput(-48,20)(0,-4){3}{\line(0,-1){4}}
\multiput(-48,20)(0,-4){4}{\line(1,0){4}}
\multiput(-48,20)(0,-4){3}{\line(0,-1){4}}
\multiput(-44,20)(0,-4){3}{\line(0,-1){4}}
\multiput(-44,20)(0,-4){4}{\line(1,0){4}}
\multiput(-40,20)(0,-4){3}{\line(0,-1){4}}
\multiput(-40,20)(0,-4){3}{\line(0,-1){4}}
\multiput(-50.8,7.7)(4,4){3}{*}
\put(-36,14){\vector(1,0){6}}
\put(-35,16){\small{$Q_4$}}
\put(-63,16){\scriptsize{$n=9$}}
\put(-63,12){\scriptsize{$\alpha_3=3$}}
\multiput(-24,20)(0,-4){5}{\line(1,0){4}}
\multiput(-24,20)(0,-4){4}{\line(0,-1){4}}
\multiput(-20,20)(0,-4){4}{\line(0,-1){4}}
\multiput(-20,20)(0,-4){4}{\line(1,0){4}}
\multiput(-20,20)(0,-4){3}{\line(0,-1){4}}
\multiput(-16,20)(0,-4){3}{\line(0,-1){4}}
\multiput(-16,20)(0,-4){4}{\line(1,0){4}}
\multiput(-12,20)(0,-4){3}{\line(0,-1){4}}
\multiput(-12,20)(0,-4){3}{\line(0,-1){4}}
\multiput(-12,20)(0,-4){2}{\line(1,0){4}}
\multiput(-12,20)(0,-4){1}{\line(0,-1){4}}
\multiput(-8,20)(0,-4){1}{\line(0,-1){4}}
\multiput(-22.8,7.7)(4,4){3}{*}
\put(-23.5,4.8){$+$}
\put(-15.5,8.8){$-$}
\put(-11.5,16.8){$+$}
\put(-6,14){\vector(1,0){6}}
\put(-5,16){\small{$Q_5$}}
\multiput(4,20)(0,-4){6}{\line(1,0){4}}
\multiput(4,20)(0,-4){5}{\line(0,-1){4}}
\multiput(8,20)(0,-4){5}{\line(0,-1){4}}
\multiput(8,20)(0,-4){4}{\line(1,0){4}}
\multiput(8,20)(0,-4){3}{\line(0,-1){4}}
\multiput(12,20)(0,-4){3}{\line(0,-1){4}}
\multiput(12,20)(0,-4){3}{\line(1,0){4}}
\multiput(12,20)(0,-4){2}{\line(0,-1){4}}
\multiput(16,20)(0,-4){2}{\line(0,-1){4}}
\multiput(16,20)(0,-4){3}{\line(1,0){4}}
\multiput(16,20)(0,-4){2}{\line(0,-1){4}}
\multiput(20,20)(0,-4){2}{\line(0,-1){4}}
\multiput(5.2,7.7)(4,4){3}{*}
\put(4.5,0.8){$+$}
\put(8.5,8.8){$-$}
\put(16.5,12.8){$+$}
\put(24,14){\vector(1,0){6}}
\put(25,16){\small{$Q_6$}}
\multiput(34,20)(0,-4){7}{\line(1,0){4}}
\multiput(34,20)(0,-4){6}{\line(0,-1){4}}
\multiput(38,20)(0,-4){6}{\line(0,-1){4}}
\multiput(38,20)(0,-4){3}{\line(1,0){4}}
\multiput(38,20)(0,-4){2}{\line(0,-1){4}}
\multiput(42,20)(0,-4){2}{\line(0,-1){4}}
\multiput(42,20)(0,-4){3}{\line(1,0){4}}
\multiput(42,20)(0,-4){2}{\line(0,-1){4}}
\multiput(46,20)(0,-4){2}{\line(0,-1){4}}
\multiput(46,20)(0,-4){3}{\line(1,0){4}}
\multiput(46,20)(0,-4){2}{\line(0,-1){4}}
\multiput(50,20)(0,-4){2}{\line(0,-1){4}}
\multiput(39.2,11.7)(4,4){2}{*}
\put(35.2,3.7){*}
\put(34.5,-3.2){$+$}
\put(53,14){\vector(1,0){6}}
\multiput(64,20)(0,-4){7}{\line(1,0){4}}
\multiput(64,20)(0,-4){6}{\line(0,-1){4}}
\multiput(68,20)(0,-4){6}{\line(0,-1){4}}
\multiput(68,20)(0,-4){3}{\line(1,0){4}}
\multiput(68,20)(0,-4){2}{\line(0,-1){4}}
\multiput(72,20)(0,-4){2}{\line(0,-1){4}}
\multiput(72,20)(0,-4){3}{\line(1,0){4}}
\multiput(72,20)(0,-4){2}{\line(0,-1){4}}
\multiput(76,20)(0,-4){2}{\line(0,-1){4}}
\multiput(76,20)(0,-4){3}{\line(1,0){4}}
\multiput(76,20)(0,-4){2}{\line(0,-1){4}}
\multiput(80,20)(0,-4){2}{\line(0,-1){4}}
\multiput(69.2,11.7)(4,4){2}{*}
\put(65.2,3.7){*}
\put(82,16){\scriptsize{$n=12$}}
\put(82,12){\scriptsize{$\alpha_3=3$}}
\end{picture}
\vskip 8mm
\caption{The construction of a partition of $12$ with three $3$-hooks.}\label{figure12}
\end{center}
\end{figure}

 Figure \ref{figure12} gives an illustration of the construction of a partition of $12$
 with three $3$-hooks from the nearly $3$-triangular partition with three parts.

Combining Corollary \ref{cor-r} and Theorem \ref{thm2}, we    obtain Theorem \ref{cor}.
We conclude this paper with a proof of Theorem \ref{thm3}.

\noindent {\it{Proof of Theorem \ref{thm3}.}} Multiplying both sides
by $(1-q)$, (\ref{generatebnk}) can be rewritten as
\begin{equation}\label{transbnk}
(1-q)\sum_{n\geq 0}b(n,k)q^n=\sum_{t\geq 1}^{\infty}q^{{t \choose
2}k^2}\frac{1-q^{tk^2}}{1-q^{tk}}-1.
\end{equation}
It is easily checked that
\begin{align}
(1-q)\sum_{n\geq 0}b(n,k)q^n=\sum_{n\geq 0}(b(n+1,k)-b(n,k))q^{n+1}.
\label{difbnk}
\end{align}

To compute  $b(n+1,k)-b(n,k)$, we denote the interval
  $[t(m,k), t(m+1,k)-1]$ by $I_m^{(k)}$. By Theorem
\ref{cor},  we see that $b(n,k)$ is determined by the interval containing $n$.
There are two cases:

\noindent Case 1: $n$ and $n+1$ belong to
the same interval $I_m^{(k)}$.
 By Theorem \ref{cor}, we have $b(n+1,k)=b(n,k)=m$, that is, $b(n+1,k)-b(n,k)=0$.

\noindent Case 2: $n$ and $n+1$ lie in two consecutive intervals $I_m^{(k)}$ and $I_{m+1}^{(k)}$,
that is, $n=t(m+1,k)-1$. By Theorem \ref{cor}, we obtain that $b(n,k)=m$ and $b(n+1,k)=m+1$.
It follows that $b(n+1,k)-b(n,k)=1$.

Combining the above two cases, we deduce that
\begin{equation*}
\sum_{n\geq 0}(b(n+1,k)-b(n,k))q^{n+1}=\sum_{m\geq 0}q^{t(m+1,k)}=\sum_{m\geq 0}q^{t(m,k)}-1.
\end{equation*}
It follows that
\begin{equation}\label{bnk-}
(1-q)\sum_{n\geq 0}b(n,k)q^n=\sum_{m\geq 0}q^{t(m,k)}-1.
\end{equation}
Recall that
\[ t(m,k)=m\left(\left\lfloor \frac{m}{k}\right\rfloor
+1\right)k-{\lfloor\frac{m}{k}\rfloor +1\choose 2}k^2.\]
 Write $m=jk+r$, where $0\leq r\leq k-1$. Then we get
 \begin{align}
 t(jk+r,k)&=(jk+r)\left(j
+1\right)k-{j +1\choose 2}k^2\notag\\
&={j+1 \choose 2}k^2+r(j+1)k.\label{tmk-}
\end{align}
Substituting \eqref{tmk-} into \eqref{bnk-}, we find that
\begin{align*}
\sum_{m\geq 0}q^{t(m,k)}-1
=&\sum_{j\geq 0}\sum_{r=0}^{k-1}q^{{j+1 \choose 2}k^2+r(j+1)k}-1\\
=&\sum_{j\geq 0}q^{{j+1 \choose 2}k^2}\frac{1-q^{(j+1)k^2}}{1-q^{(j+1)k}}-1.
\end{align*}
In view of \eqref{bnk-}, we obtain that
\begin{equation*}
(1-q)\sum_{n\geq 0}b(n,k)q^n
=\sum_{t\geq 1}q^{{t \choose 2}k^2}\frac{1-q^{tk^2}}{1-q^{tk}}-1.
\end{equation*}
This completes our proof.
\qed

Using the Jacobi triple product identity \cite[Eq. $1.6.1$]{Gasper}, we may express the
generating function $\sum_{n\geq 0}b(n,k)q^n$ in the following form:
\begin{equation*}
\frac{1}{1-q}\left(\frac{(q^{2k^2};q^{2k^2})_\infty}{(q^{k^2};q^{2k^2})_\infty}+
\frac{1}{2}\sum_{r=1}^{k-1}(-q^{rk};q^{k^2})_\infty(-q^{k^2-rk};q^{k^2})_\infty
(q^{k^2};q^{k^2})_\infty-\frac{k+1}{2}\right).
\end{equation*}

\vspace{0.5cm}
 \noindent{\bf Acknowledgments.} This work was supported by  the 973
Project, the PCSIRT Project of the Ministry of Education,  and the
National Science Foundation of China.

\end{document}